# Optimal condition of boundary flex control for the systems governed by Boussinesq equation with the press boundary condition and mixed boundary condition

Gol Kim [a]

[a] Center of Natural Sciences, University of Sciences, Kwahakdong-1, Unjong District, Pyongyang, DPR Korea (E-mail:golkim124@yahoo.com)

**Abstract.** In this paper, the boundary flex control problem of non stationary equation governing the coupled mass and heat flow of a viscous incompressible fluid in a generalized Boussinesq approximation by assuming that viscosity and heat conductivity are dependent on temperature has been studied. The boundary condition for velocity of fluid is non -standard boundary condition: specifically the case where dynamical pressure is given on some part of the boundary and the boundary condition for temperature of fluid is mixed boundary condition has been considered.
The optimal condition has been derived. Then, Pontryagin's maximum principle in the special case has been derived.
**Keyword:** Boussinesq equation, boundary flex control, press boundary condition, mixed boundary condition

## 1. Introduction

The extremal problem for Navier-Stokes equation and optimal control of fluid dynamical equation are studied by various authors (for example; [6-10], [17-19]).

In [22] the existence of time optimal controls for the Boussinesq equation has been obtained and derived Pontryagin's maximum principle of time optimal control problem governed by the Boussinesq equation

In [23] an optimal control problem governed by a system of nonlinear partial differential equations modeling viscous incompressible flows submitted to variations of temperature has been consider. A generalized Boussinesq approximation has been used. The existence of the optimal control as well as first order optimality conditions of Pontryagin type by using the Dubovitskii-Milyutin formalism has been obtained.

In [24] the stationary Boussinesq equations describing the heat transfer in the viscous heat-conducting fluid under inhomogeneous Dirichlet boundary conditions for velocity and mixed boundary conditions for temperature are considered. The optimal control problems for these equations with tracking-type functionals are formulated. A local stability of the concrete control problem solutions with respect to some disturbances of both cost functionals and state equation is proved.

In [25] the boundary control problems of the model of heat and mass transfer in a viscous incompressible heat conducting fluid has been considered. The model consists of the Navier-Stokes equations and the convection-diffusion equations for the substance concentration and the temperature that are nonlinearly related via buoyancy in the Oberbeck–Boussinesq approximation and via convective mass and heat transfer.

In [25] control problems for stationary magnetohydrodynamic equations of a viscous heat-conducting fluid under mixed boundary conditions has been considered.

In [15, 16] the Karhunen-Loeve Galerkin method for the inverse problems of Boussinesq equation have been studied.

In [20] the problem of stabilization of the Boussinesq equation via internal feedback controls has been studied. In [21] the problem of local internal controllability of the Boussinesq system has been studied.

In this paper, the boundary flex optimal control for the evolution equation governing the coupled mass and heat flow of a viscous incompressible fluid in a generalized Boussinesq approximation by assuming that viscosity and heat conductivity are dependent on temperature has been studied. The boundary condition for velocity of fluid is non -standard boundary condition: specifically the case where dynamical pressure is given on some part of the boundary will be considered. The boundary condition for temperature of fluid is mixed boundary condition.

The existence of the optimal control has been proved. Then the optimal condition has been derived.

Let $\Omega \subset R^N$ (N=2, 3) be a bounded domain with smooth boundary $\Gamma$. Let $\Gamma$ be divided by into two parts $\Gamma_1, \Gamma_2$ such that $\Gamma = \Gamma_1 \cup \Gamma_2 (\Gamma_1 \cap \Gamma_2 = \emptyset, \Gamma_2 \neq \emptyset)$ T. ($0<T<\infty$) is given number.

We denote $Q = \Omega \times (0,T)$, $\Sigma_i = \Gamma_i \times (0,T)$ (i = 1,2), $\Sigma = \Gamma \times (0,T)$.

We assume that the state of control systems is given by non- stationary Boussinesq equation with the dynamical pressure condition and mixed boundary condition on some part of the boundary as follows:

$$\begin{cases} \dfrac{\partial z}{\partial t} - v\Delta z + (z,\nabla)z + \beta \xi\, w = -grad\ \pi; & Q \\ div z = 0; & Q \\ \dfrac{\partial w}{\partial t} - k\Delta w + (z,\nabla)w = 0; & Q \\ z_\varsigma = 0,\ \pi + \dfrac{1}{2}|z|^2 = v_1,\ w = 0; & \Sigma_1 \\ z = 0,\ -k\dfrac{\partial w}{\partial n} = v_2; & \Sigma_2 \\ z(0) = z_0, w(0) = w_0; & \Omega \end{cases} \quad (1.1)$$

where $\Omega \subset R^N$ (N=2, 3) is a bounded domain with smooth boundary $\Gamma$. $\Gamma$ is divided by into two parts $\Gamma_1, \Gamma_2$ such that $\Gamma = \Gamma_1 \cup \Gamma_2 (\Gamma_1 \cap \Gamma_2 \neq \Phi)$, T ($0<T<\infty$) is given number.

We denote $Q = \Omega \times (0,T)$, $\Sigma_i = \Gamma_i \times (0,T)$ (i = 1,2), $\Sigma = \Gamma \times (0,T)$ and $n$ note the outer normal vector to $\Gamma$. In the Eqs. (1)-(6) $z(x,t) \in R^N$ denotes the velocity of the fluid at point $x \in \Omega$ at time $t \in [0,T]$; $\pi(x,t) \in R$ is the hydrostatic pressure; $w(x,t) \in R$ is temperature; g is the gravitational vector, and $v > 0$ and $k > 0$ are kinematic viscosity and thermal conductivity, respectively; $\beta$ is a positive constant associated to the coefficient of volume expansion; $v_1$ and $v_2$ are the given functions on $\Sigma_1$ and $\Sigma_2$ respectively. In Eqe.(1) $\beta > 0$ is the coefficient of volume expansion and $\xi$ is the gravitational function.

The expressions $\nabla, \Delta$ and div denote the gradient, Laplacian and divergence operators, respectively (sometimes, we will also denote the gradient operator by grad); $i$ th component in Cartesian coordinates of $(z,\nabla)z$ is given by

$$((z,\nabla)z)_i = \sum_{j=1}^{N} z_j \frac{\partial z_i}{\partial x_j}, \text{ also } (z,\nabla)w = \sum_{j=1}^{N} z_j \frac{\partial w}{\partial x_j}$$

In the boundary condition (4) $z_\varsigma = z - z_n$, $z_n = (z \cdot n)n$.

There are the results of research of Boussinesq equation and the generalized Boussinesq system with nonlinear thermal diffusion in [1-4, 23]. But boundary conditions of those papers are homogeneous.

We assume that the cost functional J[v] is given as following:

$$J[v] = J[v,y] = N_1 \int_0^T \int_{\Gamma_1} r_1(s,t) z_n(s,t) ds dt + N_2 \int_{\Gamma_2} r_2(s,t) \frac{\partial w}{\partial n} ds dt$$

where $N_1, N_2 > 0$ are given real numbers and $r_1(x,t) \in [L^2(\Sigma_1)]^N$, $r_2(x,t) \in L^2(\Sigma_2)$ are given functions .And $v = \{v_1, v_2\}$, $y = \{z, w\}$.

$\int_{\Gamma_1} r_1(s,t) z_n(s,t) ds = \int_{\Gamma_1} (r_1(s,t), z_n(s,t)) ds$ and $(r_1(s,t), z_n(s,t))$ denote the scalar product in $[L^2(\Sigma_1)]^N$.

Then, the problem that we are going to considered is to find the $v_* \in U_a$ satisfying:

$$\inf_{v \in U_a} J[v] = J[v_*] \qquad (1.2)$$

Here, we assume that the admissible control sets $U_\alpha = U_{1\alpha} \times U_{2\alpha}$ are defined such as:

$$U_{1\alpha} = \{v_1 \mid v_1 \in L^2(0,T;(L^2(\Gamma_1))^N), 0 < \alpha_1(x,t) \leq v_1(x,t) \leq \beta_1 (almost\ everywhere)\} \quad (1.3)$$

$$U_{2\alpha} = \{v_2 \mid v_2 \in L^2(0,T;L^2((\Gamma_2)), 0 < \alpha_2 \leq v_2(x,t) \leq \beta_2 (almost\ everywhere)\} \qquad (1.4)$$

$\alpha_i(x,t), \beta_i(x,t)(i=1,2)$ are given functions in the function space $L^2(\Sigma_i)$.
$v_1(x,t) = (v_{11}(x,t), v_{12}(x,t), \cdots, v_{1N}(x,t))$ and expression $0 < \alpha_1(x,t) \leq v_1(x,t) \leq \beta_1$ means that $0 < \alpha_1(x,t) \leq v_{1i}(x,t) \leq \beta_1 (\forall i \in N)$

We denote $v = \{v_1, v_2\}$ and $y(v) = \{z(v), w(v)\}$.

The established optimization problem (1.1)-(1.4) is an optimal boundary flex control problem.

For the convenient, we have assumed that control parameters are the flex pressure $v_1$ on the boundary $\Sigma_1$ and the heat flex $v_2$ on the boundary $\Sigma_2$

To illustrate the example of extremal condition (1.2), we can take functions $r_1(x,t)$ and $r_2(x,t)$ as following;

$$r_1(x,t) = \begin{cases} I & ; x \in \overline{\Gamma}_1 \subset \Gamma_1 \\ 0 & ; x \notin \overline{\Gamma}_1 \subset \Gamma_1 \end{cases}, \qquad r_2(x,t) = \begin{cases} 1 & ; x \in \overline{\Gamma}_2 \subset \Gamma_2 \\ 0 & ; x \notin \overline{\Gamma}_2 \subset \Gamma_2 \end{cases}$$

Where, I is a unit vector. Then, the optimization problem (1-2) is described the problem that fluid flex passed the boundary $\overline{\Gamma}_1$ under the restriction for the flex pressure and heat flex passed the boundary $\overline{\Gamma}_2$ under the restriction for the heat flex must minimize.

## 2. Preliminaries

In this article the functions are either R or $R^N$ ($N=2$ or $N=3$) and as usual simplification, sometimes we will not distinguish them in our notations; the difference will be clear from the context. The $L^2(\Omega)$-product and norm are denoted by $(\cdot,\cdot)$ and $|\cdot|$ respectively: the $H^m(\Omega)$ norm is denoted by $\|\cdot\|_m$. Here $H^m(\Omega) = W^{m,2}(\Omega)$ is the usual Sobolev spaces (see [1] for their properties); $H^{-1}(\Omega)$ denotes the dual spaces of $H_0^1(\Omega)$. $H^0(\Omega)$ is the same as $L^2(\Omega)$ and $\|\cdot\|_0$ is the same as $L^2$-norm $|\cdot|$. $D(0,T)$ is the class of $C^\infty$-functions with compact support in $(0,T)$. $D(0,T)'$ are its associated spaces of distribution.

If B is a Banach space, we denoted by $L^q(0,T;B)$ the Banish space of the B-valued functions defined in the interval $(0,T)$ that are $L^q$-integrable.

Now we introduce some spaces such as;

$$\begin{cases} D = \{\phi \mid \phi \in (c^\infty(\Omega))^N, div\phi(x) = 0 (x \in \Omega), \phi_\varsigma(x) = 0 (x \in \Gamma_2)\} \\ H = \text{completion of D under the } [L^2(\Omega)]^N\text{-norm} \\ V = \text{completion of D under the } [H^1(\Omega)]^N - \text{norm} \\ D_{\Gamma_1} = \{\phi : \varphi \in C^\infty(\Omega), \varphi(x) = 0 (x \in \Gamma_1)\} \\ \widetilde{H} = \text{closure of } D_{\Gamma_1} \text{ in } L^2(\Omega) \\ W = \text{closure of } D_{\Gamma_1} \text{ in } H^1(\Omega) \end{cases} \qquad (2.1)$$

Naturally, the norm of H or $\tilde{H}$ is also denoted by $|\cdot|$, and the norm of $V$ or $W$ is denoted by $\|\cdot\|$ as well. The dual product between $V^*$ and $V$ or $W^*$ and $W$ (also the inner product in $H^{-1}(\Omega)$ and $H_0^1(\Omega)$) are denoted by $<\cdot,\cdot>$.

Now, generally, we assume that $v_1 \in L^2(0,T;(H^{-1/2}(\Gamma_1))^N)$, $v_2 \in L^2(0,T;H^{-1/2}(\Gamma_2))$. Then, we shall obtain the expression of operator form the weak solution for state equation (1.1) Suppose that $\{z,w\}$ is a classical solution of (1.1). Multiplier the first equation of (1.1) by $\psi \in V$, integrate by parts over $\Omega$ and take the boundary condition into account to get

$$\frac{d}{dt}(z,\psi) + v\, a_1(z,\psi) + b(z,z,\psi) + (\beta \xi w,\psi) = <v_1, \psi_n>_{\Gamma_1}, \quad \psi_n = (\psi \cdot n)n. \tag{2.2}$$

Multiplier the second equation of (1.1) by $\varphi \in W$, integrate by parts over $\Omega$, and take the boundary conditions into account to obtain

$$\frac{d}{dt}(w,\varphi) + k a_2(w,\varphi) + c(z,w,\varphi) = <v_2,\varphi>. \tag{2.3}$$

Now let $\chi \in C^1[0,T]$ be a function such that $\chi(T) = 0$. Multiplier Equations (2.2) and (2.3) by _ $\chi$ respectively, and integrate by parts to yield

$$\int_0^T [-(z(t),\psi\chi'(t))dt + \int_0^T [va_1(z(t),\psi\chi(t)) + b(z(t),z(t),\psi\chi(t)) + \beta(w(t)\xi,\psi\chi(t))]dt$$

$$= \int_0^T \int_{\Gamma_1} (v_1(t),\psi_n \chi(t))dsdt + (z_0,\psi)\chi(0) \tag{2.4}$$

$$\int_0^T [-(w(t),\varphi\chi'(t))dt + \int_0^T [ka_2(w(t),\varphi\chi(t)) + c(z(t),w(t),\varphi\chi(t))]dt$$

$$= \int_0^T \int_{\Gamma_{21}} (v_2(t),\varphi\chi(t))dsdt + (w_0,\varphi)\chi(0) \tag{2.5}$$

where $z(\cdot,t) = z(t)$, $w(\cdot,t) = w(t)$, $v_i(\cdot,t) = v_i(t)$, $i = 1,2$ by abuse of notation without confusion from the context, and

$$(z,\psi) = \sum_{j=1}^N \int_\Omega z_j(x,\cdot)\psi_j(x)dx, \quad (w,\varphi) = \int_\Omega w(x,\cdot)\varphi(x)dx$$

$$a_1(u,\psi) = (rotz(x,\cdot), rot\psi) = \int_\Omega (\nabla \times u)\cdot(\nabla \times \phi)dx \quad b(z,z,\phi) = \int_\Omega (rotz(x,\cdot) \times z(x,\cdot))\psi(x)dx,$$

$$a_2(w,\varphi) = \sum_{j=1}^N \int_\Omega \frac{\partial w(x,\cdot)}{\partial x_j}\cdot \frac{\partial \varphi}{\partial x_j}dx, \quad c(z,w,\varphi) = \sum_{j=1}^N \int_\Omega z_j(x,\cdot)\frac{\partial w(x,\cdot)}{\partial x_j}\varphi(x)dx$$

The following Lemmas 2.1-2.3 can be obtained by the Sobolev inequalities and the compactness theorem. We can also refer to theorem 1.1 of [31] on page 107 and lemmas 1.2, 1.3 of [31] on page 109 (see also chapter 2 of [13]). The similar arguments can be also found in lemmas 1, 5 of [30].

**Lemma 2.1.** The bilinear forms $a_1(\cdot,\cdot)$ and $a_2(\cdot,\cdot)$ are coercive over V and W respectively. That is, there exist constants $c_1, c_1' > 0$ such that

$$a_1(z,z) \geq c_1\|z\|^2, \forall z \in V \text{ and } a_2(w,w) \geq c_1'\|z\|^2, \forall \in W$$

**Lemma 2.2.** The trilinear forms $b(\cdot,\cdot,\cdot)$ is a linear continuous functional on $[H^1(\Omega)]^N$. That is, there exist a constant c2 > 0 such that

$$|b(u,v,w)| \leq c_2\|u\|\cdot\|v\|\cdot\|w\|, \quad \forall u,v,w \in [H^1(\Omega)]^N$$

Moreover, the following properties hold true

Moreover, the following properties hold true

(i) $b(u, v, v) = 0, \forall u, v \in V$

(ii) $b(u, v, w) = -b(u, w, v), \quad \forall u \in V, \forall v, w \in [H^1(\Omega)]^N$

(iii) If $u_m \to u$ weakly on V and $v_m \to v$ strongly on H, then $b(u_m, v_m, w) \to b(u, v, w), \forall u, v \in V \; \forall w \in V$

**Lemma 2.3.** The tri-linear form $c(\cdot,\cdot,\cdot)$ is a linear continuous functional defined on $V \times W \times W$. That is, there exist a constant $c_3 > 0$ as following:

$$|c(z, w, \varphi)| \leq c_3 \|z\| \cdot \|w\| \cdot \|\varphi\|, \quad \forall z \in V, \quad \forall w, \varphi \in W$$

Moreover, the following properties hold true

(i) $c(z, w, w) = 0, \forall z \in V, \forall w \in W$

(ii) $c(z, w, \varphi) = -c(z, \varphi, w), \quad \forall z \in V, \forall w, \varphi \in W$

(iii). If $z_m \to z$ weakly on V and $w_m \to w$ strongly on $\tilde{H}$, then $c(z_m, w_m, \varphi) \to c(z, w, \varphi), \forall z \in V, w \in \tilde{H}, \varphi \in W$.

**Definition 1.** Let $Y \equiv Z \times W = (L^2(0,T:V) \cap L^\infty(0,T:H)) \times (L^2(0,T:W) \cap L^\infty(0,T:\tilde{H}))$.

Suppose that $v_1 \in (L^2(0,T;(H^{-1/2}(\Gamma_1))^N)$, $v_2 \in (L^2(0,T;H^{-1/2}(\Gamma_2)))$, $z_0 \in H$, $w_0 \in \tilde{H}$, $g \in L^\infty(\Omega)$.

The pair $y = \{z, w\}$ is said to be a weak solution of (1.1) if it satisfies

$$\begin{cases} y = \{z, w\} \in Y, z' \in L^1(0,T,V^*), w' \in L^1(0,T:W^*) \\ (z',\psi) + \nu a_1(z,\psi) + b(z,z,\psi) + (\beta\xi w,\psi) = <v_1,\psi_n>_{\Gamma_1}, \forall \psi \in V \\ (w',\varphi) + ka_2(w,\varphi) + c(z,w,\varphi) = <v_2,\varphi>_{\Gamma_2}, \forall \varphi \in W \\ z(0) = z_0, w(0) = w_0 \end{cases} \quad (2.7)$$

Next, we reformulate Equation (2.7) into the operator equation. To this purpose, it is noticed that for a fixed $\psi \in V$, the functional $\psi(\in V) \to a_1(z,\psi)$ is linear continuous. So there exists an $A_1 z \in V^*$ such that

$$<A_1 z, \psi> = a_1(z,\psi), \quad \forall \psi \in V \quad (2.8)$$

Similarly, for fixed $u, v \in V$, $w \in V \to b(u, v, w)$ is a linear continuous functional on V. Hence there exist a $B(u, v) \in V^*$ such that

$$<B(u, v), w> = b(u, v, w), \quad \forall w \in V \quad (2.9)$$

We denote $B(u) = B(u, u)$. Define

$$L_1(\psi) = (v_1, \psi_n)_{\Gamma_1} = \int_{\Gamma_1} v_1 \psi_n ds, \quad \forall \psi \in V. \quad (2.10)$$

Then, fixed $v_1 \in (L^2(0,T;(H^{-1/2}(\Gamma_1))^N)$, the functional $\psi(\in V) \to L_1(\psi) = (v_1,\psi_n)_{\Gamma_1}$ is linear continuous. So that there exist constant $c_4 > 0$ such that

$$\|L_1\psi\| < c_4 \|\psi\|_V, \forall \psi \in V.$$

So there exists an $H_1 v_1 \in V^*$ such that

$$<H_1 v_1, \psi> = (v_1, \psi_n)_{\Gamma_1}, \forall \psi \in V \quad (2.11)$$

With these operators at hand, we can write the second equation of (2.7) as

$$\frac{dz}{dt} + \nu A_1 z + B(z) + \beta\xi w = H_1 v_1 \quad (2.12)$$

Similarly, we have

$$<A_2 w, \varphi> = a_2(w,\varphi), \quad <C(z,w),\varphi> = c(z,w,\varphi), \quad A_2 w, C(z,w) \in W^* \quad (2.13)$$

Define

$$L_2(\varphi) = (v_2, \varphi)_{\Gamma_2} = \int_{\Gamma_2} v_2 \varphi ds, v_2 \in L^2(\Gamma_2), \forall \varphi \in W \tag{2.14}$$

Then, fixed $v_2 \in (L^2(0,T; H^{-1/2}(\Gamma_2))$ the functional $\varphi(\in W) \to L_2(\varphi) = (v_2,\varphi)_{\Gamma_{21}}$ is linear continuous.

Then, the operator $L_2$ is a linear continuous functional defined on $W$ and so there exists constant $c_5 > 0$ such that

$$\|L_2 \varphi\| \le c_5 \|\varphi\|_W, \forall \varphi \in W.$$

So there exists an $H_2 v_2 \in W^*$ such that

$$< H_2 v_2, \varphi > = (v_2, \varphi)_{\Gamma_2}, \forall \varphi \in V \tag{2.15}$$

By these operators defined above, we can write the first equation of (2.7) as

$$\frac{dw}{dt} + kA_2 w + C(z,w) = H_2 v_2 \tag{2.16}$$

Combining (2.12) and (2.16), we can write (2.7) in the abstract evolution equation as follows:

$$\begin{cases} \frac{dz}{dt} + vA_1 z + B(z) + \beta \xi w = H_1 v_1, \\ \frac{dw}{dt} + kA_2 w + C(z,w) = H v_2, \\ z(0) = z_0, w(0) = w_0 \end{cases} \tag{2.17}$$

**Lemma 2.4.** If $z \in L^2(0,T;V)$, then $B(z) \in L^1(0,T;V^*)$; and if $w \in L^2(0,T;W)$, then $C(z,w) \in L^1(0,T;W^*)$.

Proof. By applying Höler inequality and compactness of embedding $H^1(\Omega) \subset L_4(\Omega)$, we obtain

$$|(B(z), \psi)| = |b(z,z,\psi)| \le c_6' \|z\| \cdot \|z\|_{L_4} \cdot \|\ \|_{L_4} \le c_6 \|z\|^2 \|\psi\|,$$

for some constants $c_6', c_6 > 0$ and hence $\|B(z)\|_{V^*} \le c_6 \|z\|^2$ which shows that $B(z) \in L^1(0,T;V^*)$. The proof is complete.
Similarly, we have

$$|(C(z,w), \varphi)| = |-c(z,\varphi,w)| \le c_7 \|z\|_{L_4} \cdot \|\varphi\| \cdot \|w\|_{L_4} \le c_8 \|z\| \cdot \|\varphi\| \cdot \|w\| \le c_9 (\|z\|^2 + \|w\|^2) \|\varphi\|$$

for some constants $c_8, c_9 > 0$ and hence $\|C(z,w)\| \le c_9 (\|z\|^2 + \|w\|^2)$ which shows that $C(z,w) \in L^1(0,T;W^*)$ for all $\varphi \in W$. The proof is complete.

We specify the constants $c_i, i = 1,2,\cdots$ used in this section in the remaining part of the paper. The following Lemma 2.5 comes from theorem 2.2 of [12] on page 220

**Definition 2.** Suppose that $z_0 \in H$, $w_0 \in \widetilde{H}$, $g \in L^\infty(\Omega)$.
The pair $\{y, v\} = \{(z,w), (v_1, v_2)\}$ is said to be admissible pair; "state-control" of extreme value problem (1.2) if it satisfies

$$\begin{cases} y = \{z, w\} \in Y, z' \in L^1(0,T,V^*), w' \in L^1(0,T:W^*) \\ \frac{dz}{dt} + vA_1 z + B(z) + \beta \xi w = H_1 v_1, \\ \frac{dw}{dt} + kA_2 w + C(z,w) = H v_2, \\ z(0) = z_0, w(0) = w_0 \\ v = \{v_1, v_2\} \in U_a \equiv U_{1a} \times U_{2a} \end{cases} \tag{2.18}$$

We denote admissible pair sets by M that is;

$$M \equiv \{(y,v) \mid y \equiv \{z,w\} \in Y = Z \times W, v = \{v_1, v_2\} \in U_a, (y,v) \text{ satisfy (2.18)}\}.$$

## 3. Optimal condition

In this section, we derive the first order necessary condition of optimal control problem (1.2). When deriving, the optimal condition in the control problem of nonlinear system, usually additional regularity condition for admissible pair is demanded. Here we assume that regularity condition for optimal state $\{z_*, w_*\}$ in the case three-dimensional domain hold true such that;

$$z_* \in L^6(0,T;V), \quad w_* \in L^6(0,T;W) \tag{3.1}$$

We assume that this regularity condition is satisfied here and in what follows.
In order to derive the optimal condition of optimal control problem, for arbitrary $\varepsilon > 0$ we introduce $\varepsilon$-approximation control problem;

**Problem $P_\varepsilon$.**

$$\inf_{v \in U_a} J[v] \tag{3.2}$$

State constraints are given such as;

$$g'_\tau + v A_1 g_\tau + B(z_*, g_\tau) + B(g_\tau, z_*) + \varepsilon B(g_\tau, g_\tau) + \beta \xi \eta_\tau = H_1(v_1 - v_{1*}) \tag{3.3}$$

$$\eta'_\tau + k A_2 \eta_\tau + C(z_*, \eta_\upsilon) + C(g_\tau, w_*) + \varepsilon C(\eta_\tau, \eta_\tau) = H_2(v_2 - v_{2*}) \tag{3.4}$$

$$g_\tau(0) = 0, \, g_\tau \in L^2(0,T;V), \, g'_\tau \in L^1(0,T;V^*) \tag{3.5}$$

$$\eta_\tau(0) = 0, \, \eta_\tau \in L^2(0,T;W), \, \eta'_\tau \in L^1(0,T;W^*) \tag{3.6}$$

Here, $\{v_{1*}, v_{2*}\}$ is an optimal control and $\{z_*, w_*\}$ is optimal state.
And, we denote $g_\tau = g(x,\tau)$, $\eta_\tau = \eta(x,\tau)$, $g'_\tau = \dfrac{\partial g(x,\tau)}{\partial \tau}$, $\eta'_\tau = \dfrac{\partial \eta(x,\tau)}{\partial \tau}$, $\tau \in [0,T]$.

**Lemma 3.1.** We assume that the condition

$$\frac{\beta \|\xi\|_\infty (\beta \|\xi\|_\infty + 1)}{v c_1} \leq \frac{k c'_1}{2}$$

is hold. Then, there exists a solution of the equation (3.3)-(3.6).
(Proof) We assume that $\{g_m(t)\}, \{\eta_m(t)\}$ is solution for the Galerkin system of equation (3.3), (3.4) that is;

$$g'_m + v A_1 g_m + Q^1_m(B(z_*, g_m) + B(g_m, z_*)) + \varepsilon B(g_m, g_m) + \beta \xi \eta_m = Q^1_m H_1(v_1 - v_{1*}), g_m(0) = 0 \tag{3.7}$$

$$\eta'_m + k A_2 \eta_m + Q^2_m(C(z_*, \eta_m) + C(g_m, w_*)) + \varepsilon C(\eta_m, \eta_m) = Q^2_m H_2(v_2 - v_{2*}), \eta_m(0) = 0 \tag{3.8}$$

Here $Q^1_m$ is projection operator from $H$ to $H^1_m$ and $H^1_m$ is the subspace generated by first m numbers "basis" element of space $H$. And $Q^2_m$ is projection operator from $\tilde{H}$ to $H^2_m$ and $H^2_m$ is the subspace generated by first $m$ numbers "basis" element of space $\tilde{H}$
By applying scalar product (3.7) by $g_m$, we obtain

$$\frac{1}{2}\frac{d}{d}|g_m|^2 + v a_1(g_m, g_m) + \langle B(g_m, z_*), g_m \rangle + (\beta \xi \eta_m, g_m) = \langle H_1(v_1 - v_{1*}), g_m \rangle$$

From here, for any given $\varepsilon > 0$, we obtain we can get that

$$\frac{1}{2}\frac{d}{d}|g_m|^2 + v c_1 a_1 \|g_m\|^2 \leq -\langle B(g_m, z_*), g_m \rangle - (\beta \xi \eta_m, g_m) + \langle H_1(v_1 - v_{1*}), g_m \rangle$$

$$\leq -\langle B(g_m, z_*), g_m \rangle + \frac{1}{2}\beta \|\xi\|_\infty (\frac{1}{\varepsilon^2}\|\eta_m\|^2 + \varepsilon^2 \|g_m\|^2) + \frac{1}{2}[\frac{1}{\varepsilon^2}|v_1 - v_{1*}|_{L^2(\Sigma_1)} + \varepsilon^2 \|g_m\|^2]$$

$$\leq \left|-\langle B(g_m, z_*), g_m \rangle\right| + \frac{1}{2}(\beta\|\xi\|_\infty + 1)\varepsilon^2 \|g_m\|^2 + \frac{1}{2}\beta\|\xi\|_\infty \frac{1}{\varepsilon^2}\|\eta_m\|^2 + \frac{1}{2}\frac{1}{\varepsilon^2}|v_1 - v_{1*}|_{L^2(\Sigma_1)}$$

From here, we obtain that

$$\frac{1}{2}\frac{d}{dt}|g_m|^2 + [vc_1 - \frac{1}{2}(\beta\|\xi\|_\infty + 1)\varepsilon^2]\|g_m\|^2 \leq \left|-\langle B(g_m, z_*), g_m \rangle\right| + \frac{1}{2}\beta\|\xi\|_\infty \frac{1}{\varepsilon^2}\|\eta_m\|^2 + \frac{1}{2}\frac{1}{\varepsilon^2}|v_1 - v_{1*}|_{L^2(\Sigma_1)}$$

Now, we put $\varepsilon^2 = \dfrac{vc_1}{(\beta\|\xi\|_\infty + 1)}$, then obtain that

$$\frac{1}{2}\frac{d}{dt}|g_m|^2 + \frac{1}{2}vc_1\|g_m\|^2 \leq c_4|v_1 - v_{1*}|_{L^2(\Sigma_1)} + \left|-\langle B(g_m, z_*), g_m \rangle\right| + \frac{1}{2}\beta\|\xi\|_\infty \frac{(\beta\|\xi\|_\infty + 1)}{vc_1}\|\eta_m\|^2$$

Here, $c_4 = \dfrac{1}{2}\dfrac{(\beta\|\xi\|_\infty + 1)}{vc_1}$.

From here, we obtain that

$$|g_m(t)|^2 + vc_1 \int_0^t \|g_m(\tau)\|^2 d\tau \leq 2c_4|v_1 - v_{1*}|^2_{L^2(\Sigma_1)} + c_{10}\int_0^t \|g_m(\tau)\|^{1+\frac{N}{4}} \|z_*\|^{\frac{N}{4}} (|z_*|\|g_m\|)^{1-\frac{N}{4}} d\tau +$$

$$+ \beta\|\xi\|_\infty \frac{(\beta\|g\|_\infty + 1)}{vc_1} \int_0^t \|\eta_m\|^2 d\tau$$

Multiplying (3.8) by $\eta_m$

$$\frac{1}{2}\frac{d}{dt}|\eta_m|^2 + ka_2(\eta_m, \eta_m) + \langle C(g_m, w_*), \eta_m \rangle = \langle H_2(v_2 - v_{2*}), \eta_m \rangle$$

Same as above, we obtain that with $c_5 = \dfrac{1}{2kc'_1}$

$$|\eta_m(t)|^2 + kc'_1 \int_0^t \|\eta_m(\tau)\|^2 d\tau \leq 2c_5|v_2 - v_{2*}|^2_{L^2(\Sigma_2)} + c_{11}\int_0^t \|g_m(\tau)\| \|w_*\|^{\frac{N}{4}} (|w_*|\|\eta_m\|)^{1-\frac{N}{4}} \|\eta_m\|^{\frac{N}{4}} d\tau$$

Here we have applied inequality as following;

$$|\phi|_{L^4} \leq k\|\phi\|^{\frac{N}{4}} |\phi|^{1-\frac{N}{4}}, \quad \forall \phi \in V( \text{ or } \phi \in W)$$

$k > 0$ is constant.

By adding above two inequalities, we obtain

$$|g_m(t)|^2 + |\eta_m(t)|^2 + vc_1 \int_0^t \|g_m(\tau)\|^2 d\tau + kc'_1 \int_0^t \|\eta_m(\tau)\|^2 d\tau \leq 2c_4|v_1 - v_{1*}|^2_{\Sigma_1} + c_5|v_2 - v_{2*}|^2_{\Sigma_2} +$$

$$+ c_{10}\int_0^t \|g_m(\tau)\|^{1+\frac{N}{4}} \|z_*\|^{\frac{N}{4}} (|z_*|\|g_m\|)^{1-\frac{N}{4}} d\tau + c_{11}\int_0^t \|g_m(\tau)\| \|w_*\|^{\frac{N}{4}} (|w_*|\|\eta_m\|)^{1-\frac{N}{4}} \cdot \|\eta_m\|^{\frac{N}{4}} d\tau$$

$$+ \beta\|\xi\|_\infty \frac{(\beta\|\xi\|_\infty + 1)}{vc_1} \int_0^t \|\eta_m\|^2 d\tau \tag{3.9}$$

Here and in what follows, we denote $|\cdot|^2_{L^2(\Sigma_1)}$ by $|\cdot|^2_{\Sigma_1}$ and $|\cdot|^2_{[L^2(\Sigma_1)]^N}$ by $|\cdot|^2_{\Sigma_2}$ simply.

Now, let estimate the terms of right-hand side of (3.9).

We apply Young's inequality and $|z_*| \leq c$, $|w_*| \leq c$ in the third term of right-hand side of (3.9). That is, in Young's inequality $ab \leq \dfrac{a^p}{p} + \dfrac{b^q}{q}$, $p > 1, \dfrac{1}{p} + \dfrac{1}{q} = 1$ if we put $p = \dfrac{8}{4-N}$, then $q = \dfrac{8}{4+N}$.

Therefore, we obtain

$$\int_0^t \|g_m\|^{1+\frac{N}{4}} \|w_*\|^{\frac{N}{4}} (w_* \|\eta_m\|)^{1-\frac{N}{4}} \|\eta_m\|^{\frac{N}{4}} d\tau = \int_0^t (\|z_*\|^{\frac{N}{4}} |g_m|^{1-\frac{N}{4}})(|z_*|^{1-\frac{N}{4}} \|g_m\|^{1+\frac{N}{4}}) d\tau \le$$

$$\le \frac{4-N}{8} \int_0^t (\|z_*\|^{\frac{N}{4} \cdot \frac{8}{4-N}} |g_m|^{\frac{4-N}{4} \cdot \frac{8}{4-N}}) d\tau + \frac{4-N}{8} \int_0^t |z_*|^{\frac{4-N}{4} \cdot \frac{8}{4+N}} \|g_m(\tau)\|^{\frac{4-N}{4} \cdot \frac{8}{4+N}} d\tau$$

$$= \frac{4-N}{8} \int_0^t (\|z_*\|^{\frac{2N}{4-N}} |g_m|^2) d\tau + \frac{4+N}{8} \int_0^t |z_*|^{\frac{2(4-N)}{4+N}} \|g_m\|^2 d\tau$$

$$\le \frac{4-N}{8} \int_0^t (\|z_*\|^{\frac{2N}{4-N}} |g_m|^2) d\tau + \frac{4+N}{8} [\frac{4-N}{4+N} \int_0^t |z_*|^2 d\tau + \frac{2N}{4+N} \int_0^t \|g_m\|^{2 \cdot \frac{4+N}{2N}} d\tau]$$

$$\le \frac{4-N}{8} \int_0^t (\|z_*\|^{\frac{2N}{4-N}} |g_m|^2) d\tau + \frac{4-N}{8} \int_0^t |z_*|^2 d\tau + \frac{N}{4} \int_0^t \|g_m(\tau)\|^{\frac{4+N}{N}} d\tau$$

$$\le \frac{4-N}{8} \int_0^t (\|z_*\|^{\frac{2N}{4-N}} |g_m|^2) d\tau + \frac{4-N}{8} \int_0^t |z_*|^2 d\tau + \nu c_1 \int_0^t \|g_m\|^2 d\tau + C(N,\nu,c_1) \int_0^t d\tau$$

$$\le \frac{4-N}{8} \int_0^t (\|z_*\|^{\frac{2N}{4-N}} |g_m|^2) d\tau + \frac{4-N}{8} \int_0^t |z_*|^2 d\tau + \nu c_1 \int_0^t \|g_m\|^2 d\tau + C(N,\nu,T,c_1) \quad (3.10)$$

Here $C(N,\nu,T)$ is constant related $N$ and $\nu, T$.

Similarly, by estimating the forth term of right-hand side of (3.9) we obtain

$$\int_0^t \|g_m\| \|w_*\|^{\frac{N}{4}} (w_* \|\eta_m\|)^{1-\frac{N}{4}} \|\eta_m\|^{\frac{N}{4}} d\tau = \int_0^t (\|w_*\|^{\frac{N}{4}} |\eta_m|^{1-\frac{N}{4}}) \cdot (|w_*|^{1-\frac{N}{4}} \|g_m\| \|\eta_m\|^{\frac{N}{4}}) d\tau \le$$

$$\le \frac{4-N}{8} \int_0^t (\|w_*\|^{\frac{N}{4} \cdot \frac{8}{4-N}} |\eta_m|^{\frac{4-N}{4} \cdot \frac{8}{4-N}}) d\tau + \frac{4+N}{8} \int_0^t |w_*|^{\frac{4-N}{4} \cdot \frac{8}{4+N}} \|g_m\|^{\frac{8}{4+N}} \|\eta_m\|^{\frac{N}{4} \cdot \frac{8}{4+N}} d\tau$$

$$= \frac{4-N}{8} \int_0^t (\|w_*\|^{\frac{2N}{4-N}} |\eta_m|^2) d\tau + \frac{4+N}{8} \int_0^t |w_*|^{\frac{2(4-N)}{4+N}} \|g_m\|^{\frac{8}{4+N}} \|\eta_m\|^{\frac{2N}{4+N}} d\tau$$

$$\le \frac{4-N}{8} \int_0^t (\|w_*\|^{\frac{2N}{4-N}} |\eta_m|^2) d\tau + \nu c_1 \int_0^t \|g_m\|^2 d\tau + C(\nu,N,c_1) \frac{4+N}{8} \int_0^t |w_*|^{\frac{2(4-N)}{4+N} \cdot \frac{4+N}{N}} \|\eta_m\|^{\frac{2N}{4+N} \cdot \frac{4+N}{N}} d\tau$$

$$= \frac{4-N}{8} \int_0^t (\|w_*\|^{\frac{2N}{4-N}} |\eta_m|^2) d\tau + \nu c_1 \int_0^t \|g_m\|^2 d\tau + C(\nu,N,c_1) \frac{N}{8} \int_0^t |w_*|^{\frac{2(4-N)}{N}} \|\eta_m\|^2 d\tau$$

$$\le \frac{4-N}{8} \int_0^t (\|w_*\|^{\frac{2N}{4-N}} |\eta_m|^2) d\tau + \nu c_1 \int_0^t \|g_m\|^2 d\tau + C(\nu,N,c_1) \frac{N}{8} \frac{4-N}{N} \int_0^t |w_*|^2 d\tau + \frac{N}{8} C(\nu,N) \int_0^t \|\eta_m\|^{2 \cdot \frac{N}{2N-4}} d\tau$$

$$= \frac{4-N}{8} \int_0^t (\|w_*\|^{\frac{2N}{4-N}} |\eta_m|^2) d\tau + \nu c_1 \int_0^t \|g_m\|^2 d\tau + \frac{4-N}{8} C(\nu,N,c_1) \int_0^t |w_*|^2 d\tau + \frac{2N-4}{8} C(\nu,N) \int_0^t \|\eta_m\|^{\frac{2N}{2N-4}} d\tau$$

$$\le \frac{4-N}{8} \int_0^t (\|w_*\|^{\frac{2N}{4-N}} |\eta_m|^2) d\tau + \nu c_1 \int_0^t \|g_m\|^2 d\tau + \frac{4-N}{8} C(\nu,N,c_1) \int_0^t |w_*|^2 d\tau +$$

$$+ kc_1' \int_0^t \|\eta_m\|^2 d\tau + \frac{2N-4}{8} \frac{N}{2N-4} C(\nu,N,c_1) \int_0^t d\tau = \frac{4-N}{8} \int_0^t (\|w_*\|^{\frac{2N}{4-N}} |\eta_m|^2) d\tau + \nu c_1 \int_0^t \|g_m\|^2 d\tau +$$

$$+ kc_1' \int_0^t \|\eta_m\|^2 d\tau + \frac{4-N}{8} C(\nu, N, c_1) \int_0^t |w_*|^2 d\tau + C(\nu, \kappa, N, T, c_1, c_1') \qquad (3.11)$$

Here $C(\nu, \kappa, N, T, c_1, c_1')$ is constant related $N$ and $\nu, k, T, c_1, c_1'$.

By applying (3.10), (3.11) in (3.9) we obtain

$$|g_m(t)|^2 + |\eta_m(t)|^2 + \nu c_1 \int_0^t \|g_m(\tau)\|^2 d\tau + kc_1' \int_0^t \|\eta_m(\tau)\|^2 d\tau \le 2c_4 |v_1 - v_{1*}|_{\Sigma_1}^2 +$$

$$2c_5 |v_2 - v_{2*}|_{\Sigma_2}^2 + c_{10} \int_0^t |z_*|^2 d\tau + c_{11} \int_0^t |w_*|^2 + C(N, \nu, T, C(\nu, \kappa, N, T, c_1)) + C(\nu, \kappa, N, T, c_1, c_1'), c_1') +$$

$$+ c_{10} \int_0^t [\|z_*\|^{\frac{2N}{4-N}} + \frac{1}{2}|\xi|^2] \cdot |g_m|^2 d\tau + c_{11} \int_0^t [\|w_*\|^{\frac{2N}{4-N}} + \frac{1}{2}] \cdot |\eta_m|^2 d\tau + \beta \|\xi\|_\infty \frac{(\beta \|\xi\|_\infty + 1)}{\nu c_1} \int_0^t \|\eta_m\|^2 d\tau$$

Now, applying the condition of the theorem, we obtain such as

$$|g_m(t)|^2 + |\eta_m(t)|^2 + \nu c_1 \int_0^t \|g_m(\tau)\|^2 d\tau + \frac{kc_1'}{2} \int_0^t \|\eta_m(\tau)\|^2 d\tau \le 2c_4 |v_1 - v_{1*}|_{\Sigma_1}^2 +$$

$$2c_5 |v_2 - v_{2*}|_{\Sigma_2}^2 + c_{10} \int_0^t |z_*|^2 d\tau + c_{11} \int_0^t |w_*|^2 + C(N, \nu, T, C(\nu, \kappa, N, T, c_1)) + C(\nu, \kappa, N, T, c_1, c_1'), c_1') +$$

$$+ c_{10} \int_0^t [\|z_*\|^{\frac{2N}{4-N}} + \frac{1}{2}|\xi|^2] \cdot |g_m|^2 d\tau + c_{11} \int_0^t [\|w_*\|^{\frac{2N}{4-N}} + \frac{1}{2}] \cdot |\eta_m|^2 d\tau \qquad (3.12)$$

Therefore,

$$|g_m(t)|^2 + |\eta_m(t)|^2 + \nu c_1 \int_0^t \|g_m(\tau)\|^2 d\tau + kc_1' \int_0^t \|\eta_m(\tau)\|^2 d\tau \le c_{12} +$$

$$c_{10} \int_0^t M_1(\tau) \cdot |g_m|^2 d\tau + c_{12} \int_0^t M_2(\tau) \cdot |\eta_m|^2 d\tau \qquad (3.13)$$

Here, $c_{12} = 2c_4 |v_1 - v_{1*}|_{\Sigma_1}^2 + 2c_5 |v_2 - v_{2*}|_{\Sigma_2}^2 + c_{10} \int_0^t |z_*|^2 d\tau + c_{11} \int_0^t |w_*|^2 + C(N, \nu, T) + C(\nu, \kappa, N, T)$ and here and in what follows, we denote

$$M_1(\tau) = \|z_*(\tau)\|^{\frac{2N}{4-N}} + \frac{1}{2}|\xi|^2, M_2(\tau) = \|w_*(\tau)\|^{\frac{2N}{4-N}} + \frac{1}{2} \qquad (3.14)$$

and denote $|\xi|_{R^N}^2$ by $|\xi|^2$ simply.

From (3.13), we can obtain

$$|g_m(t)|^2 + |\eta_m(t)|^2 + \nu c_1 \int_0^t \|g_m(\tau)\|^2 d\tau + kc_1' \int_0^t \|\eta_m(\tau)\|^2 d\tau \le c_{12} + c_{10} \int_0^t M_1(\tau) \cdot |g_m|^2 d\tau +$$

$$c_{11} \int_0^t M_2(\tau) \cdot |\eta_m|^2 d\tau \le c_{12} + c_{13} \int_0^t [M_1(\tau) \cdot |g_m|^2 + M_2(\tau) \cdot |\eta_m|^2] d\tau$$

Here, $c_{13} = \max\{c_{10}, c_{11}\}$.

Then, we can find a nonnegative function $M_3(\tau)$ that is bounded and integrable in a.e. $\tau \in [0, T]$ such that;

$$\max\{M_1(\tau), M_2(\tau)\} \le M_3(\tau), \text{in a.e. } \tau \in [0, T].$$

Therefore, from above inequality, we obtain

$$|g_m(t)|^2 + |\eta_m(t)|^2 + \nu c_1 \int_0^t \|g_m(\tau)\|^2 d\tau + kc_1' \int_0^t \|\eta_m(\tau)\|^2 d\tau \leq c_{12} + c_{13} \int_0^t M_3(\tau)[|g_m|^2 + |\eta_m|^2] dt \quad (3.15)$$

By (3.15), we obtain $|g_m(t)|^2 + |\eta_m(t)|^2 \leq c_{14} + c_{14} \int_0^t M_2(\tau) \cdot [|g_m|^2 + |\eta_m|^2] d\tau$

By Gronwall inequality we can obtain

$$|g_m(t)|^2 + |\eta_m(t)|^2 \leq c_{14} \exp\{c_{14} \int_0^t M_2(\tau) d\tau\}$$

From here, we obtain

$$\|g_m\|_{L^\infty(0,T;H)} \leq \text{const}, \quad \|\eta_m\|_{L^\infty(0,T;\tilde{H})} \leq \text{const} \quad (3.16)$$

Again, replacing $t \to T$ in (3.15), we obtain

$$\|g_m\|_{L^2(0,T;V)} \leq \text{const} \quad \|\eta_m\|_{L^2(0,T;W)} \leq \text{const} \quad (3.17)$$

Accordingly, from (3.7), (3.8) we obtain

$$\|g_m'\|_{L^1(0,T;V^*)} \leq \text{const}, \quad \|\eta_m'\|_{L^1(0,T;W^*)} \leq \text{const} \quad (3.18)$$

Therefore, from (3.16)- (3.18) we can conclude that there exists a subsequence of $\{g_m\}, \{\eta_m\}$ which we denote by the same symbols, such that

$$g_m \to g_\varepsilon, \ast\text{-weakly in } L^\infty(0,T:H), g_m \to g_\varepsilon \text{ weakly in } L^2(0,T:V)$$
$$g_m \to g_\varepsilon, \text{ strongly in } L^2(0,T:H) \quad (3.19)$$
$$\eta_m \to \eta_\varepsilon \ast\text{-weakly in } L^\infty(0,T:\tilde{H}), \eta_m \to \eta_\varepsilon \text{ weakly in } L^2(0,T:W)$$
$$\eta_m \to \eta_\varepsilon \text{ strongly in } L^2(0,T:\tilde{H}) \quad (3.20)$$

And

$$g_\varepsilon \in L^\infty(0,T;H) \cap L^2(0,T;V), g_\varepsilon' \in L^1(0,T;V^*),$$
$$\eta_\varepsilon \in L^\infty(0,T;\tilde{H}) \cap L^2(0,T;W), g_\varepsilon' \in L^1(0,T;W^*).$$

From (3.19), (3.20), by taking the limit as $m \to \infty$ in (36), (37), we obtain that $\{g_\varepsilon\}, \{\eta_\varepsilon\}$ is the solution of (3.3)-(3.6). □

**Lemma 3.2.** There is at the least of solutions $(g_\varepsilon, \eta_\varepsilon; v_{1\varepsilon}, v_{2\varepsilon})$ of problem $P_\varepsilon$

The proof of this lemma is just same as lemma 1

Now, let derive the optimal condition of approximation control problem $P_\varepsilon$.

Let examine $\{(v_{1*} + \varepsilon(v_1 - v_{1*}), v_{2*} + \varepsilon(v_2 - v_{2*})), (z_* + \varepsilon g_*, w_* + \varepsilon \eta)\}$ is admissible pair of problem (1.2).

First of all, $\{(v_{1*}, v_{2*}), (z_*, w_*)\}$ being admissible pair, we obtain

$$z_*' + \nu A_1 z_* + B(z_*, z_*) + \beta \xi w_* = H_1 v_{1*} \quad (3.21)$$
$$w_*' + k A_2 w_* + C(w_*, z_*) = H_2 v_{2*} \quad (3.22)$$

Now, multiplying (3.3) by $\varepsilon$ and adding (3.21), we obtain

$$(z_* + \varepsilon g_\varepsilon)' + \nu A_1(z_* + \varepsilon g_\varepsilon) + B(z_* + \varepsilon g_\varepsilon, z_* + \varepsilon g_\varepsilon) + \beta \xi(w_* + \varepsilon \eta_\varepsilon) = H_1(v_{1*} + \varepsilon(v_1 - v_{1*})) \quad (3.23)$$

Multiplying (3.4) by $\varepsilon$ and adding (3.22), arranging, we obtain

$$(w_* + \varepsilon \mu_\varepsilon)' + k A_2(w_* + \varepsilon \mu_\varepsilon) + C(z_* + \varepsilon g_\varepsilon, w_* + \varepsilon \mu_\varepsilon) = H_2(v_{2*} + \varepsilon(v_2 - v_{2*})) \quad (3.24)$$

From (3.23), (3.24), we conclude that $\{(v_{1*} + \varepsilon(v_1 - v_{1*}), v_{2*} + \varepsilon(v_2 - v_{2*})), (z_* + \varepsilon g_*, w_* + \varepsilon \eta)\}$ is admissible pair.

Therefore, we can conclude such that

$$I(z_* + \varepsilon g_\varepsilon, w_* + \varepsilon \eta_\varepsilon) - I(z_*, w_*) \geq 0$$

Accordingly,

$$N_1\int_0^T\int_{\Gamma_1}r_1(x,t)(z_*+\varepsilon g_\varepsilon)_n \,ds\,dt + N_2\int_0^T\int_{\Gamma_2}r_2(x,t)\frac{\partial}{\partial n}(w_*+\varepsilon\eta_\varepsilon)\,ds\,dt -$$

$$-N_1\int_0^T\int_{\Gamma_1}r_1(x,t)(z_*)_n\,ds\,dt - N_2\int_0^T\int_{\Gamma_2}r_2(x,t)\frac{\partial w_*}{\partial n}\,ds\,dt =$$

$$=\varepsilon N_1\int_0^T\int_{\Gamma_1}r_1(x,t)(g_\varepsilon)_n\,ds\,dt + \varepsilon\int_{\Gamma_1}r_1(x,t)(g_\varepsilon n)ds\,dt + \varepsilon N_2\int_0^T\int_{\Gamma_2}r_2(x,t)\frac{\partial \eta_\varepsilon}{\partial n}\,ds\,dt \geq 0$$

From here, we obtain the optimal condition of approximation control problem $P_\varepsilon$ such as;

$$N_1\int_0^T\int_{\Gamma_1}r_1(x,t)(g_\varepsilon)_n\,ds\,dt + N_2\int_0^T\int_{\Gamma_2}r_2(x,t)\frac{\partial \eta_\varepsilon}{\partial n}\,ds\,dt \geq 0 \tag{3.25}$$

**Lemma 3.3.** We assume that $\{(v_{1*},v_{2*}),(z_*,w_*)\}$ is optimal pair. Then there exists the solution $\{g,\eta\}$ of equations;

$$g' + vA_1 g + B(z_*,g) + B(g,z_*) + \beta\xi\eta = H_1(v_1 - v_{1*}), \quad g(0) = 0 \tag{3.26}$$

$$\eta g + kA_2\eta + C(z_*,\eta) + B(\eta,w_*) = H_2(v_2 - v_{2*}), \quad \eta(0) = 0 \tag{3.27}$$

$$g \in L^2(0,T;V), \quad g' \in L^2(0,T;V^*) \tag{3.28}$$

$$\eta \in L^2(0,T;W), \quad \eta' \in L^2(0,T;W^*) \tag{3.29}$$

and inequalities

$$N_1\int_0^T\langle H_1 r_1, g\rangle dt + N_2\int_0^T\langle H_2 r_2, \eta\rangle dt \geq 0 \tag{3.30}$$

are satisfied. Here, we denote $g = g(x,t)$, $\eta = \eta(x,t)$, $g' = \frac{\partial g(x,t)}{\partial t}$, $\eta' = \frac{\partial \eta(x,t)}{\partial t}$, $\tau \in [0,T]$.

Proof. By applying scalar product (3.3) by $g_\varepsilon$, we have

$$\frac{1}{2}\frac{d}{d}|g_\varepsilon|^2 + va_1(g_\tau,g_\tau) + \langle B(g_\varepsilon,z_*),g_\varepsilon\rangle + (\beta\xi\eta_\varepsilon,g_\varepsilon) = \langle H_1(v_1-v_{1*}),g_\varepsilon\rangle$$

By applying scalar product (3.4) by $\eta_\varepsilon$, we have

$$\frac{1}{2}\frac{d}{d}|\eta_\varepsilon|^2 + ka_2(\eta_\tau,\eta_\tau) + \langle C(g_\varepsilon,w_*),\eta_\varepsilon\rangle = \langle H_2(v_2-v_{2*}),\eta_\varepsilon\rangle$$

Repeating the process of proof of lemma 3.1, we obtain the estimations same as (3.16)- (3.18) such that;

$$\|g_\varepsilon\|_{L^\infty(0,T;H)} \leq const, \quad \|g_\varepsilon\|_{L^2(0,T;V)} \leq const \tag{3.31}$$

$$\|\eta_\varepsilon\|_{L^\infty(0,T;\tilde{H})} \leq const, \quad \|\eta_\varepsilon\|_{L^2(0,T;W)} \leq const \tag{3.32}$$

$$\|g'_k\|_{L^1(0,T;V^*)} \leq const, \quad \|\eta'_k\|_{L^1(0,T;W^*)} \leq const \tag{3.33}$$

Therefore, from (3.31)-(3.33) we can conclude that there exists a subsequence of $\{g_\varepsilon\}$, $\{\eta_\varepsilon\}$ which we denote by the same symbols, such that

$$g_m \to g_\varepsilon, *\text{-weakly in } L^\infty(0,T:H), g_m \to g_\varepsilon \text{ weakly in } L^2(0,T:V)$$

$$g_m \to g_\varepsilon, \text{ strongly in } L^2(0,T:H) \tag{3.34}$$

$$\eta_m \to \eta_\varepsilon *\text{-weakly in } L^\infty(0,T:\tilde{H}), \eta_m \to \eta_\varepsilon \text{ weakly in } L^2(0,T:W)$$

$$\eta_m \to \eta_\varepsilon \text{ Strongly in } L^2(0,T:\tilde{H}) \tag{3.35}$$

And

$$g_\varepsilon \in L^\infty(0,T;H)\cap L^2(0,T;V), g'_\varepsilon \in L^1(0,T;V^*),$$

$$\eta_\varepsilon \in L^\infty(0,T;\tilde{H})\cap L^2(0,T;W), g'_\varepsilon \in L^1(0,T;W^*).$$

By taking the limit as $m \to \infty$ in (3.3), (3.4), from (3.34), (3.35) we obtain that $\{g, \eta\}$ is the solution of (3.26), (3.27). By considering (3.25), we obtain (3.30) □

Now, let introduce conjugate systems such as;

$$-p' + \nu A_1 p + B(z_*, p) + B(p, z_*) + C(q, w_*) = H_1 r_1, \quad p(T) = 0 \tag{3.36}$$

$$-q' + k A_2 q + C(z_*, q) + \beta \xi\, p = H_2 r_2, \quad q(T) = 0 \tag{3.37}$$

$$p \in L^2(0, T; V), p' \in L^2(0, T; V^*), q \in L^2(0, T; W), q' \in L^2(0, T; W^*) \tag{3.38}$$

**Lemma 3.4.** There exists solution of conjugate systems (3.36)-(3.38).

Proof. The existence of solution for linear problem (3.36)-(3.38) is proved by taking the limit as $m \to \infty$ in Galerkin system of equation

$$-p'_m + \nu A_1 p_m + Q^1_m(B(z_*, p_m) + B(p_m, z_*)) + Q^2_m C(q_m, w_*) = Q^1_m H_1 r_1, \quad p_m(T) = 0$$

$$-q'_m + k A_2 q_m + Q^2_m C(z_*, q_m) + \beta \xi\, p_m = Q^2_m H_2 r_2, \quad q_m(T) = 0$$

Here $Q^1_m$ is projection operator from $H$ to $H^1_m$ and $H^1_m$ is the subspace generated by first $m$ numbers "basis" element of space $H$. And $Q^2_m$ is projection operator from $\tilde{H}$ to $H^2_m$ and $H^2_m$ is the subspace generated by first $m$ numbers "basis" element of space $\tilde{H}$. □

**Theorem.** We assume that $(\{v_{1*}, v_{2*}\}, \{z_*, w_*\})$ is arbitrary optimal pair. Then there exists $\{p, q\}$ which satisfy the equations and inequalities such as;

$$z'_* + \nu A_1 z_* + B(z_*, z_*) + \beta \xi w_* = H_1 v_{1*}, \quad z(0) = z_0 \tag{3.39}$$

$$w'_* + k A_2 w_* + C(w_*, z_*) = H_2 v_{2*}, \quad w(0) = w_0 \tag{3.40}$$

$$z_* \in L^2(0, T; V), w_* \in L^2(0, T; W) \tag{3.41}$$

$$-p' + \nu A_1 p + B(z_*, p) + B(p, z_*) + C(q, w_*) = H_1 r_1, \quad p(T) = 0 \tag{3.42}$$

$$-q' + k A_2 q + C(z_*, q) + \beta \xi\, p = H_2 r_2, \quad q(T) = 0 \tag{3.43}$$

$$p \in L^2(0, T; V), q \in L^2(0, T; W) \tag{3.44}$$

$$N_1 \int_0^T \int_{\Gamma_1} p_n (v_1 - v_{1*}) ds dt + N_2 \int_0^T \int_{\Gamma_1} q(v_2 - v_{2*}) ds dt \leq 0, \tag{3.45}$$

$$\forall v_1(x, t) \in [\alpha_1(x, t), \beta_1(x, t)], \forall v_2(x, t) \in [\alpha_2(x, t), \beta_2(x, t)]$$

(Proof) Multiplying (3.36) by $g(t)$ and integrating from 0 to T for $t$, we have

$$\int_0^T \langle -p' + \nu A_1 p + B(z_*, p) + B(p, z_*) + C(q, w_*), g(t) \rangle dt = \int_0^T \langle H_1 r_1, g(t) \rangle dt$$

By applying the integration by parts in right hand of this expression, we obtain

$$\int_0^T \langle g' + \nu A_1 g + B(z_*, g) + B(g, z_*), p \rangle dt + \int_0^T \langle C(q, w_*), g(t) \rangle dt = \int_0^T \langle H_1 r_1, g(t) \rangle dt \tag{3.46}$$

Multiplying (3.37) by $\eta(t)$ and integrating from 0 to T for $t$, we have

$$\int_0^T \langle -q' + k A_2 q + C(z_*, q) + \beta \xi\, p, \eta(t) \rangle dt = \int_0^T \langle H_2 r_2, \eta(t) \rangle dt$$

By applying the integration by parts in right hand of this expression, we obtain

$$\int_0^T \langle \eta' + k A_2 \eta + C(z_*, \eta) + B(g, z_*), q \rangle dt + \int_0^T \langle \beta \xi\, \eta, p \rangle dt = \int_0^T \langle H_2 r_2, \eta(t) \rangle dt \tag{3.47}$$

By adding (3.46) and (3.47), we obtain

$$\int_0^T \langle g' + \nu A_1 g + B(z_*, g) + B(g, z_*) + \beta \xi \eta, p \rangle dt + \int_0^T \langle \eta' + k A_2 \eta + C(z_*, \eta) + C(g, w_*), q(t) \rangle dt =$$

$$= \int_0^T \langle H_1 r_1, g(t)\rangle dt + \int_0^T \langle H_2 r_2, \eta(t)\rangle dt$$

By applying (3.25)- (3.30) in the right hand of above expression, we have

$$\int_0^T \langle -H_1(v_1 - v_{1*}), p(t)\rangle dt + \int_0^T \langle -H_2(v_2 - v_{2*}), q(t)\rangle dt \geq 0$$

From here, by applying the definition of $H_1$ and $H_2$, we obtain (3.45) □

**Corollary 3.1.** (Pontryagin's maximum principle in the special case )
We assume that control $v_2 \equiv 0$ in optimal control problem (1.2) and $(\{v_{1*}\}, \{z_*, w_*\})$ is arbitrary optimal pair. Then optimal control $v_{1*}$ satisfy maximal condition in almost everywhere of $\Sigma_1$ such as

$$v_{1*}(x,t)(p \cdot n)(x,t) = \sup_{\alpha_1(x,t) \leq \mu_1 \leq \beta_1(x,t)} (\mu_1 \cdot (p \cdot n)(x,t)) \quad , \text{a.e.}(x,t) \in \Sigma \qquad (3\text{-}48)$$

Like the preceding, we assume that control $v_2 \equiv 0$ in optimal control problem (1.2) and $(\{v_{2*}\}, \{z_*, w_*\})$ is arbitrary
optimal pair. Then optimal control $v_{1*}$ satisfy maximal condition in almost everywhere of $\Sigma_2$ such as

$$v_{2*}(x,t)q(x,t) = \sup_{\alpha_2(x,t) \leq \mu_2 \leq \beta_2(x,t)} (\mu_2 \cdot q(x,t)) \quad , \text{a.e.}(x,t) \in \Sigma \qquad (3\text{-}49)$$

Here $\{p, q\}$ is the solution of conjugate system.
(Proof) First of all, let prove (3-48).
We assume that $E$ is whole points of $\Sigma_1$ which is Lebesgue points for $p_n$, $v_{1*} \in L^2(\Sigma_1)$ and $(x_0, t_0) \in E$.
In the case $v_2 \equiv 0$, Optimal condition (3.45) is simplified such as

$$\int_0^T \int_{\Gamma_1} p_n(v_1 - v_{1*}) ds dt \leq 0, \quad \forall v_1(x,t) \in [\alpha_1(x,t), \beta_1(x,t)], \qquad (3.50)$$

Then, we put in (5.50) $v_1(x,t)$ such as;

$$v_1(x,t) = v_{1*} + \aleph_j(\mu_1 - v_{1*}(x,t)), \quad \forall \mu_1 \in [\alpha_1(x,t), \beta_1(x,t)]$$

Here, $\aleph_j$ is characteristic function in certain neighborhood of point $(x_0, t_0)$ in $\Sigma_1$ which converge to point $(x_0, t_0)$ when $j \to \infty$.
Now, by dividing (3.50) by $\int_{\Sigma_1} \aleph_j ds dt$, we obtain

$$\frac{1}{\int_{\Sigma_1} \aleph_j ds dt} \cdot \int_{\Sigma_1} \aleph_j p_n (\mu_1 - v_{1*}(x,t)) ds dt \leq 0$$

By taking the limit as $j \to \infty$ in above inequality, we obtain

$$p_n(\mu_1 - v_{1*}(x,t)) \leq 0, \forall \mu_1 \in [\alpha_1(x,t), \beta_1(x,t)] \quad, a.\ e.(x,t) \in \Sigma_1$$

From here we obtain maximal principle (3.48). By same as above method, we obtain (3.49) □
**Remark.** In fact, the function $p_n(x,t) = (p \cdot n) \in L^2(\Sigma_1)$ is switching function of fluid for optimal control.
Seeing that, by (3.48)

if $p_n(x,t) > 0$ ,then $v_{1*}(x,t) = \beta_1(x,t)$,

$$\text{if } p_n(x,t) < 0, \text{then } v_{1*}(x,t) = \alpha_1(x,t).$$

Similarly, we obtain

$$\text{if } q(x,t) > 0, \text{then } v_{2*}(x,t) = \beta_2(x,t),$$
$$\text{if } p_n(x,t) < 0, \text{then } v_{2*}(x,t) = \alpha_2(x,t).$$

Therefore, the function $q(x,t) \in L^2(\Sigma_2)$ is switching function of heat flow for optimal control.

In the case $p_n\big|_{\Sigma_1} \neq 0$, a.e. $(x,t) \in \Sigma_1$  $q\big|_{\Sigma_2} \neq 0$, a.e. $(x,t) \in \Sigma_2$, equalities (3.48),(3.49) denote "Bang-Bang Principle" which well known in the optimal control theory.

# References

[1] S. A. Lorca, J.L.Boldrini, The initial value problem for a generalized Boussinesq model, Nonlinear Analysis, 36(4), 1999, 457-480,

[2] G.Galiano, Spatial and time localization of solutions of the Boussinesq system with nonlinear thermal diffusion, Nonlinear Analysis, 42, 2000, 423-438,

[3] J. I. Diaz, G. Galiano, On the Boussinesq system with nonlinear thermal diffusion, Nonlinear Analysis, 30, 1997, 3255-3263,

[4] S. A. Lorca, J.L. Boldrini, The initial value problem for a generalized Boussinesq model: regularity and global existence of strong solutions, Matemática Contemporánea, 11, 1996, 175-201

[5] S.A. Lorca, J.L. Boldrini, Stationary solutions for generalized Boussinesq model, J. Differential Equations, 124(2), 1996, 201-225

[6] K. Óeda, On the initial value problem for the heat convection equation of Boussinesq approximation in a time-dependent domain, Proc. Japan Acad.64, Sec.A143-146, 1988

[7] H. Mormoto, Non- stationary Boussinesq equations, J. Fac. Sci. Univ. Tokyo, Sec.1A Math., 39 61-75, 1992

[8] Tang Xianjiang, On the existence and uniqueness of the solution to the Navier-Stokes Equations, Acta Mathematica Scientia, 15(3), 342-351, 1995

[9] R. Temam, Navier-Stokes Equations, Studies in Mathematics and its Applications, 2 North-Holland, Amsterdam, 1984.

[10] M.Shinbrot,W.P.Kotorynski, The initial value problem for a viscous heat-conducting fluid, J. Math. Anal.Appl.45, 1-22, 1974

[11] M.Shinbrot, W.P. Kotorynski, The initial value problem for a viscous heat-conducting fluid, J. Math. Anal. Apple. 45, 1974, 1-22,

[12] J. L. Lions, Controle des systemes distribues singuliers, Gauthier-villars, Bordas, Paris, 1983.

[13] R. Temam, Navier-Stokes Equations, Studies in Mathematics and its Applications, North-Holland, Amsterdam, 1984

[14] J. L. Lions, Quelques methods de resolution des problems aux non lineares, Dunod, paris, 1968

[15] H.M.Park, W.S.Jung, The Karhunen-Loeve Galerkin method for the inverse natural convection problems, International Journal of Heat and Mass transfer, 44, 155-167, 2001

[16] H.M.Park et al, An inverse natural convection problem of estimating the strength of a heat source, International Journal of Heat and Mass transfer, 42, 4259 - 4273, 1999.

[17] A.A.Illarionov, Asymptotics of solution to the optimal control problem for time-independent Navier-Stokes Equations, Computational Mathematics and Mathematical Physics, 41(7), 1045-1056, 2001

[18] L. Steven Hou, Thomas P. Svobodny T.P., Optimization problem for the Navier-Stokes equations with regular boundary conditions, Journal of Mathematical Analysis and Applications, 177 ( 2) 342-367, 1993.

[19] Thomas Bewley, Roger Temam, Mohammed Ziane, Existence and uniqueness of optimal

control to the Navier-Stokes equations, C. R. Acad. Sci. Paris, t.330, Serie 1, 1007-1011, 2000.

[20] Wang, Gengsheng, Stabilization of the Boussinesq equation via internal feedback controls, Nonlinear Analysis, 52(2), 2002, 485-508,

[21] Wang, Lijuan; Wang, Gengsheng, Local internal controllability of the Boussinesq system, Nonlinear Analysis, 53(5), 2003, 637-655.

[22] Shugang Li, Gengsheng Wang, The time optimal control of the Boussinesq equations Numerical Functional Analysis and Optimization, 24(1-2), 2003, 163-180

[23] Jose Luiz Boldrini, Enrique Fern'andez-Cara, Marko Antonio Rojas-Medar, An optimal control problem for a generalized Boussinesq model: The time dependent case, Rev. Mat. Complut. 20(2), 2007, 339–366

[24] Gennady Alekseev, Dmitry Tereshko, Stability of optimal controls for the stationary, International Journal of Differential Equations, Volume 2011, 2011, 1-28,

[25] Gennady Alekseev1, Dmitry Tereshko1 , Vladislav Pukhnachev, Boundary control problems for Oberbeck–Boussinesq model of heat and mass transfer, Advanced Topics in Mass Transfer, Edited by Mohamed El-Amin, 485-512, 2011

[26] G. V. Alekseev , R. V. Brizitskii, Control problems for stationary magnetohydrodynamic equations of a viscous heat-conducting fluid under mixed boundary conditions, Computational Mathematics and Mathematical Physics, 45(12), 2005, 2049–2065.

[27] R. Temam, Navier-Stokes Equations: theory and numerical analysis, AMS Chelsea Publishing, American Mathematical Society・Providence, Rhodel Island, 2001

[28] G.V. Alekseev and R.V. Brizitskii, Control problems for stationary magnetohydrodynamic equations of a viscous heat-conducting uid under mixed boundary conditions, Comput. Math. Math. Phys., 45(2005), 2049-2065.

[29] R. Temam, Navier-Stokes Equations: Theory and Numerical Analysis, AMS, Providence, Rhodel Island, 2001.

[30] G.V. Alekseev and R.V. Brizitskii, Control problems for stationary magnetohydrodynamic equations of a viscous heat-conducting uid under mixed boundary conditions, Comput. Math. Math. Phys., 45(2005), 2049-2065.

[31] R. Temam, Navier-Stokes Equations: Theory and Numerical Analysis, AMS, Providence, Rhodel Island, 2001.